\documentclass[conference]{IEEEtran}
\usepackage{amsmath}
\usepackage{algorithm}
\usepackage{algorithmicx}
\usepackage{graphicx}

\newtheorem{theorem}{Theorem}[section]

\newcommand{\qed}{\nobreak \ifvmode \relax \else
      \ifdim\lastskip<1.5em \hskip-\lastskip
      \hskip1.5em plus0em minus0.5em \fi \nobreak
      \vrule height0.75em width0.5em depth0.25em\fi}

\ifCLASSINFOpdf
\else
\fi
\begin{document}

\newcommand{\pe}{\psi}
\def\d{\delta}
\def\ds{\displaystyle}
\def\e{{\epsilon}}
\def\eb{\bar{\eta}}
\def\enorm#1{\|#1\|_2}
\def\Fp{F^\prime}
\def\Kc{{\cal K}}
\def\norm#1{\|#1\|}
\def\wb{{\bar w}}
\def\zb{{\bar z}}
\def\vulminlp{{\sc vulnerability-minlp}}
\def\theproblem{(\ref{opt:obj})--(\ref{opt:binary})}
\def\milp{{\sc milp}}
\def\minlp{{\sc minlp }}
%
\title{An Implicit Optimization Approach for Survivable Network Design}

\author{\IEEEauthorblockN{Richard Li-Yang Chen}
\IEEEauthorblockA{Quantitative Modeling and Analysis\\
Sandia National Laboratories\\
Livermore, California 94551--0969\\
Email: rlchen@sandia.gov}
\and
\IEEEauthorblockN{Amy Cohn}
\IEEEauthorblockA{Industrial and Operations Engineering\\
University of Michigan\\
Ann Arbor, Michigan 48109--2117\\
Email: amycohn@umich.edu}
\and
\IEEEauthorblockN{Ali Pinar}
\IEEEauthorblockA{Quantitative Modeling and Analysis\\
Sandia National Laboratories\\
Livermore, California 94551--0969\\
Email: apinar@sandia.gov}}


%


\maketitle

\begin{abstract}
We consider the problem of designing a network of minimum cost while satisfying a prescribed survivability criterion. The survivability criterion requires that a feasible flow must still exists (i.e. all demands can be satisfied without violating arc capacities) even after the disruption of a subset of the network's arcs. Specifically, we consider the case in which a disruption (random or malicious) can destroy a subset of the arcs, with the cost of the disruption not to exceed a disruption budget. This problem takes the form of a tri-level, two-player game, in which the network operator designs (or augments) the network, then the attacker launches a disruption that destroys a subset of arcs, and then the network operator attempts to find a feasible flow over the residual network. We first show how this can be modeled as a two-stage stochastic program from the network operator's perspective, with each of the exponential number of potential attacks considered as a disruption scenario. We then reformulate this problem, via a Benders decomposition, to consider the recourse decisions implicitly, greatly reducing the number of variables but at the expense of an exponential increase in the number of constraints. We next develop a cut-generation based algorithm. Rather than \emph{explicitly} considering each disruption scenario to identify these Benders cuts, however, we develop a bi-level program and corresponding separation algorithm that enables us to \emph{implicitly} evaluate the exponential set of disruption scenarios. Our computational results demonstrate the efficacy of this approach.
\end{abstract}


%
\IEEEpeerreviewmaketitle

\section{Introduction}
\label{sec:intro}

Society depends heavily on networked systems such as  the electric power grid, water and gas distribution systems, communication networks, and transportation systems. This reliance makes it crucial to fortify and ensure the security of these networks. The Northeast blackout in 2003 is  frequently used as  an example of the severe consequences  of an infrastructure failure. 

Many of these networks are also congested, with growth in demand exceeding expansions in capacity. This forces systems to operate close to their boundaries of feasibility. At the same time, rapid technology developments have made these systems so complex that human expertise is no longer sufficient to secure operations. Automated tools must instead be developed to detect system vulnerabilities.

Many recent research efforts have focused on vulnerability analysis for critical infrastructure networks. \cite{Wein2005} analyzed a bioterror attack on the food supply. \cite{Chen2007} examined the affects of one or more arc failures in a transportation network in terms of network travel time or generalized travel cost increase as well as the behavioral responses of users due to the failure in the network. Efforts such as these have produced significant enhancements in models and algorithms for \emph{identifying} vulnerability. The next logical step, which is the goal of this paper, is to improve systems to \emph{reduce} these vulnerabilities, abating the risk of failure.

It is often useful to represent physical networks with mathematical networks. For example, nodes can represent generators, load points, and junctions in a power system, with arcs representing transmission lines. Supplies and demands at each node correspond to generation and load, and arc capacities represent transmission limits. Mathematical constraints then approximate the physical dynamics of flow across the network. 

The underlying structure of a network can be critical in enabling us to understand and remedy its vulnerabilities. For example in \cite{Pinar2010}, it was shown that even though nonlinear equations  govern the flow of power through a network, vulnerabilities of the underlying power system can be identified by investigating the structural properties of the network. In this paper, we specifically exploit network structure to solve the problem of identifying the minimum-cost set of arc capacities to install in the network while ensuring that a feasible flow will be possible even after an arbitrary disruption of arcs that is limited only by a disruption budget $\Gamma$.

Formally, we study the following problem:

{\it 
Given a network $G=(N,E)$, a set of candidate arcs $E_c$, demand/supply on each node, capacities of all arcs, the costs of constructing each arc in $E_c$, the cost of disrupting each arc in $E\cup E_c$, and a disruption budget $\Gamma$,  find a set of arcs $E_n \subset E_c$ whose cumulative cost is  minimum, such that for any arc subset $X$ whose cumulative disruption cost is  less than $\Gamma$, a feasible flow, satisfying all demand, exist on the network $G'=(N,(E\cup E_n)\setminus X)$. 
} 

This problem poses significant challenges, largely due to the nesting of multiple optimization problems.  For a given network and disruption, the third level problem is to minimize the disruption cost.  At the second level, therefore, the attacker maximizes the impact of the disruption, subject to anticipating the third-level optimization. Finally, at the first level, the network operator must optimize the network design subject to the anticipation of the attacker's optimal response.

To formulate this problem as a traditional mixed-integer linear program requires the enumeration of all possible attacks, which may be exponentially large. For example, in the simplest case where the budget limits the number of attacked edges to be at most $k$, there will be approximately $\binom{|E \cup E_c|}{k}$ attacks to evaluate. Thus, this approach is intractable for all but very small instances, and tractability depends on the use of an alternative approach that can identify feasible solutions to the network design problem, i.e. solutions that can survive any attack within the disruption budget, without having to explicitly consider all disruption scenarios. This is the focus of our paper.
 
The rest of the paper is organized as follows. Section \ref{sec2} presents the single level network flow problem, the bilevel network disruption problem, and the trilevel network design problem.  In Section \ref{sec3}, we discuss cutting plane procedures and develop an implicit optimization approach for solving the trilevel network design problem. Section \ref{sec_comp_exp} presents our experimental design and computational results, and Section \ref{sec5} concludes the paper.

\section{Models}\label{sec2}

Consider an undirected network $G(N,E)$ with node set $N$ and \emph{undirected} arc (candidate and existing)  set $E$.  To differentiate between capacities, which are undirected, and flows, which are directed, we denote the undirected link between node $i$ and $j$ as $\{i,j\}$ and the two corresponding directed arcs as $(i,j)$ and $(j,i)$.  Let arc set $A$ represent the set of all directed arcs corresponding to the set of undirected arcs $E$, that is, $A=\{(i,j),(j,i)|\{i,j\}\in E\}$. Each arc  $\{i,j\} \in E$ is associated with a construction cost $c_{ij}$ and a capacity $u_{ij}$.  If $b_i>0$, then node $i$ is a supply node, if  $b_i<0$ then $i$ is a demand node, and  $b_i=0$ for all other (transshipment) nodes.   Without loss of generality, we assume that $\sum_{i \in N} b_i =0$.

The goal is to design a minimum cost network such that a feasible flow exists under \emph{any} disruption within the disruption budget $\Gamma$.  In the first level, the network operator determines the network design (or augmentation).  We define a binary variable $x_{ij}$,  so that
\[
x_{ij} = \left\{ \begin{array}{ll} 	1, 	& \textmd{if arc $\{i,j\}$ is constructed, }\; 	 \\
                                                   	 0, 	& {\rm otherwise,}\; 		
 \end{array}\right.  \forall \{i,j\} \in E
\]
For simplicity of notation, we do not make a distinction between existing and candidate arcs.  For an \emph{existing} arc $\{i,j\} \in E$, we set $x_{ij}$ to be one and fixed its construction cost $c_{ij} = 0$. 

The cost to the attacker of disrupting arc $\{i,j\} \in E$ is given by $r_{ij}$.  We define a binary variable $d_{ij}$,  so that
\[
d_{ij} = \left\{ \begin{array}{ll} 	1, 	& \textmd{if arc $(i,j)$ is disrupted, }\; 	 \\
                                                   	 0, 	& {\rm otherwise,}\; 		
 \end{array}\right.  \forall \{i,j\} \in E
\]

The arc disruption cost may reflect the probability of random failure or the cost of disruption to an attacker.  We assume that the amount of disruption that can be realized is constrained by a disruption budget $\Gamma$.  That is

\begin{eqnarray}
r^Td \leq \Gamma.  \nonumber
\end{eqnarray}

\subsection{Trilevel Optimization Framework}

We begin by presenting an explicit formulation of the problem as a tri-level optimization problem. Recall the three levels comprising the problem:

\begin{enumerate}
\item  In the first level, the network operator (for example, the Independent System Operator of a power system) makes decisions about which arcs to either add or augment (for example, by adding sensors, redundancy, or other protective measures to improve the impenetrability of the arc).
\item In the second level, a  disruption (with cost $\leq \Gamma$) occurs on the newly-augmented network such that a  subset of the arcs are destroyed.
\item  In the third level, the network operator responds to the disruption by seeking to find a feasible recourse, with flow able to satisfy demand without violating the capacity constraints on the surviving arcs.
\end{enumerate}

The goal of the overall problem is to find the minimum-cost first stage decisions subject to the constraint that the third stage problem will be feasible under any second stage decisions. 

To present our explicit formulation of this tri-level problem, we begin by assuming fixed values for a network design $(x)$ and disruption scenario $(d)$. The third level problem is to then minimize the  penalty $(p\gamma)$ associated with failing to fully satisfy demand. If the optimal objective value is zero, then the third level problem is in fact feasible. Note that the penalty term $p$ must be chosen to be large enough such that the penalty cost for an infeasible third stage problem will dominate the objective function when making first stage decisions, ensuring that in the final solution the third stage will be feasible for any valid disruption in the second stage.

\begin{subequations}
\begin{eqnarray}
 &\hskip -1cm  \displaystyle \min_{f,\gamma}&  p \hskip 0.02cm \gamma \label{mod_op_obj} \\
& \hskip -1cm \displaystyle \mathrm{s.t.}&   \hskip -0.2cm \sum_{j:(i,j)\in A}  f_{ij} - \sum_{j:(j,i)\in A} f_{ji} = b_i (1-\gamma)  \quad \forall i \in N \label{mod_op_flow_bal}\\
&&f_{ij}, f_{ji} \leq u_{ij}(x_{ij}-d_{ij}) \quad \forall \{i,j\} \in E \label{mod_op_flow_ub} \\
&& f_{ij} \geq 0 \quad \forall (i,j) \in A \label{mod_op_flow_lb} \\
&& \gamma \geq 0 \label{mod_op_alpha_nneg}
\end{eqnarray}
\end{subequations}

In this statement of the third level of the problem, note that $(x)$ and $(d)$ are input parameters, passed from the two upper levels, rather than decision variables. Objective (\ref{mod_op_obj}) is to minimize the fraction of demand that goes unmet (i.e. the load shed).   (\ref{mod_op_flow_bal}) are flow balance constraints that require the net flow into and out of a node to be equal to the demand (or supply). It is possible that given a network design $(x)$ there may not exist a feasible flow post disruption $(d)$.  Thus we introduce a scaling variable $\gamma \in [0,1]$ corresponding to fraction of demand shed (e.g. not satisfied).  An $\gamma$ value of zero indicates that \emph{all} of the demand is satisfied and a value of one indicates that \emph{none} of the demand is satisfied.  (\ref{mod_op_flow_ub}) and (\ref{mod_op_flow_lb}) are arc capacity upper and lower bounds, respectively. Note that the upper bound on arc $\{i,j\}$ depends on whether the arc is part of the network as defined in the first stage problem $(x_{ij})$ as well as whether the arc has been disrupted as defined in the second stage problem $(d_{ij})$.   We assume, for notational simplicity, that flow lower bounds are zero but a more general model with nonnegative lower bound is applicable as well.

We now step up one level to the second level. We assume that the network is vulnerable to disruptions due to both malicious and natural  causes, and we wish to secure the network against all possible disruptions.  We assume that when an arc is disrupted, all of its capacity is lost.  In this sense, the disruption we are considering is a worst-case disruption. For critical infrastructure protection, it is imperative to look at worst-case scenarios.  One reason for this is that critical infrastructure can easily be the target of a malicious attack by an intelligent adversary.  Another is that even highly unlikely events can cause huge disturbances, as evidenced by the Northeast blackouts \cite{Bienstock2010}\cite{Pinar2010}, due to the paramount importance of these critical infrastructures from the perspectives of both security and the economy.

Given that we want to ensure against any possible disruption (with cost $\leq \Gamma$), we assume that for any given network design (defined by the decision variables $(x)$), the most damaging disruption will realize.  We refer to this problem as the \emph{network disruption problem} (\textsc{ndp}) and formulate it as follows.

\begin{subequations}\label{ndp_bi-level}
\begin{eqnarray}
 \hspace*{-10ex}&\displaystyle \max_{d} & \min_{f ,\gamma}  \  \ \ p \hskip 0.02cm \gamma \label{mod_ndp_obj}  \\
& \mathrm{s.t.}&  \hskip -0.2cm \sum_{j:(i,j)\in A} \! f_{ij} - \!\!\sum_{j:(j,i)\in A}\! f_{ji} \! = \! b_i (1-\gamma) \;  \forall i\! \in\! N \label{mod_ndp_flow_bal}\\
& & f_{ij}, f_{ji} \leq u_{ij}(x_{ij}-d_{ij}) \quad \forall \{i,j\} \in E  \label{mod_ndp_arc_ub}\\
& & f_{ij} \geq 0 \quad \forall (i,j) \in A  \label{mod_ndp_arc_lb}\\
&& \gamma \geq 0 \\
& & \hskip -0.8cm r^T d \leq \Gamma \label{mod_ndp_budget}\\
& & \hskip -0.8cm d \in\{0,1\}^n \label{modsip_attack_decision}
\end{eqnarray}
\end{subequations}

The objective (\ref{mod_ndp_obj}) is to choose the disruption scenario for which the network operator's corresponding optimal recourse will have the highest penalty cost.  The constraints look the same as in the third stage problem, except that $(d)$ is now a decision variable rather than an input parameter. In addition, constraint (\ref{mod_ndp_budget}) is added to restrict the feasible disruption scenario to be one in which the disruption budget is not violated. The disruption cost is assumed to be an additive function of the disrupted arcs. Observe, that if $x_{ij}$ is equal to one then $d_{ij}$ can take on either one or zero.  If $x_{ij}$ is zero, then constraints (\ref{mod_ndp_arc_ub}) and (\ref{mod_ndp_arc_lb}) together force the attack decision on the corresponding component to be zero as well.  This is equivalent to the statement that non-existent components are not disruptable.  

Problem of the form (\ref{ndp_bi-level}) are commonly referred to in the literature as \emph{network interdiction problems}.  For example,  \cite{Wollmer1964} proposed an algorithm for performing sensitivity analysis on maximum flow networks.   More recently, \cite{Lim2006} examines the multicommodity network interdiction problem under discrete and continuous attacks.

Work on network interdiction problems has brought a lot of insights to the vulnerabilities of networks, and has paved the way for a higher objective: How do we build/augment networks in the first place, so as to limit their vulnerability? This new problem, which embeds the bi-level network interdiction problem as constraints, yields a tri-level optimization problem.

Tri-level optimization problems are extremely challenging and, typically, cannot be solved without much difficulty.  \cite{Brown2006} proposed general bi-level and tri-level programming models for defending critical infrastructure. \cite{Smith2007} examined the problem of augmenting a network under various disruption scenarios, assuming disruptions may be selected based on arc capacities, initial flows, and minimizing the maximum profit from transmitting flows.  \cite{Yao2007} proposed a nested bi-level programming approach for the augmentation of electric power grids.

In the network design problem which we examine, we embed the bi-level \textsc{ndp} problem, yielding a tri-level optimization problem, which we refer to as the network augmentation problem (\textsc{nap}).

\begin{subequations}
\begin{eqnarray}
 &\hskip -1.8cm \displaystyle \min_{x } &  \hskip -0.55cm \max_{d} \hskip 0.3cm  \min_{f ,\gamma}  \ \sum_{\{i,j\} \in E} c_{ij}x_{ij} + p \gamma \label{nap_obj} \\
& \hskip -1.8cm  \mathrm{s.t.}&  \hskip 0.3cm \sum_{j:(i,j)\in A} \!\! f_{ij} - \!\!\!\sum_{j:(j,i)\in A} \! \! f_{ji} \! =\! b_i (1-\gamma)  \;\;  \forall i \! \in\! N\label{nap_flow_bal}\\
&&  \hskip 0.5cm f_{ij}, f_{ji} \leq u_{ij}(x_{ij}-d_{ij}) \quad \forall \{i,j\} \in E  \label{nap_cap_ub}\\
&&  \hskip 0.5cm f_{ij} \geq 0 \quad \forall (i,j) \in A \label{nap_cap_lb} \\
&&  \hskip 0.5cm \gamma \geq 0 \\
&&  \hskip -0.4cm r^T d \leq \Gamma \label{modsip_budget} \label{nap_budget}\\
&&  \hskip -0.4cm d \in\{0,1\}^n \label{modsip_attack_decision} \label{nap_disruption_var}\\
&& \hskip -1.4cm  x \in\{0,1\}^n \nonumber
\end{eqnarray}
\end{subequations}

The objective (\ref{nap_obj}) is to minimize the total network design cost as well as the penalty associated with the worst-case disruption, given the network design $(x)$. By setting the penalty cost $p$ adequately large, the first stage decisions will always yield a network design that can survive a feasible attack, so long as such a network design exists.  Observe that the constraints of \textsc{nap} are the same as those of \textsc{ndp} except that $(x)$ is now a decision variable instead of an input parameter. Finally, note that the non-negativity constraints (\ref{nap_cap_lb}) prevent an arc from being attacked in the second stage unless that arc is created (or is pre-existing) in the first stage, else the right-hand side of (\ref{nap_cap_ub}) would be negative.

\section{Solution Approaches}\label{sec3}

We begin by making the observation that given a finite disruption budget $\Gamma$ (in the simplest case, this might be an upper bound on the number of arcs that could be disrupted), the number of possible disruption scenarios is finite (although it may be extremely large).  Let $S$ be the total number of possible disruption scenarios with cost less than or equal to the given disruption budget $\Gamma$.

One approach is to \emph{explicitly} consider all possible disruption scenarios $(s=1,\cdots,S)$ in determining the optimal network augmentation. By explicitly considering all possible disruption scenarios, we can reduce the tri-level optimization model to a single \emph{mixed integer linear program} (\textsc{milp}), although it may be extremely large.  We denote this model as the \emph{extensive form} (\textsc{ef}) and state it as follows.

\begin{subequations}\label{ef}
\begin{eqnarray}
  &\hskip -2cm \displaystyle{\min_{x,f, \gamma,\theta}}  & \hskip -0.7cm \sum_{\{i,j\} \in E} c_{ij}x_{ij} + p \hskip 0.02cm \theta \label{ef_obj}\\
&\hskip -2cm \mathrm{s.t.}  & \hskip -0.5cm   \gamma^s \leq \theta \quad \forall s \in S \label{ef_scen_op_cost}\\
&&\hskip -0.8cm \sum_{j:(i,j)\in A}  f_{ij}^s - \sum_{j:(j,i)\in A } f_{ji}^s = b_i (1-\gamma^s) \nonumber \\
&&\hskip 3.5cm  \forall i \in N, s\in S \label{ef_flow_bal}  \\
&&\hskip -0.7cm  f_{ij}^s, f_{ji}^s \leq u_{ij} (x_{ij} - d_{ij}^s) \quad \forall \{i,j\} \in E, s\in S \label{ef_arc_cap} \\
&&\hskip -0.7cm f_{ij}^s \geq 0 \quad \forall (i,j) \in A, s\in S  \label{ef_flow_nneg}\\
&&\hskip -0.7cm \theta, \gamma^s \geq 0\\
&&\hskip -0.7cm x_{ij} \in\{0,1\} \quad \forall \{i,j\} \in E \label{ef_x_integ}
\end{eqnarray}
\end{subequations}

Observe that in \textsc{ef},  $d^s$ is an $|E|$-dimensional vector (of parameters) corresponding to disruptions of each arc  $(\{i,j\} \in E)$ in a disruption scenario $s$. The objective (\ref{ef_obj}) is to minimize the network design cost and shortage cost under the worst-case disruption scenario, which will naturally vary with the network design $(x)$. Constraints (\ref{ef_scen_op_cost}) ensure that the operating cost $\theta$ is lower bounded by the operating cost under each disruption scenario.  Since the objective is to minimize total cost, $\theta$ will attain a value equal to the maximum operating cost across all scenarios for any given network design $(x)$. Constraints (\ref{ef_flow_bal}) are flow balance constraints for each node-disruption scenario pair.  Constraints (\ref{ef_arc_cap}) are arc capacity constraints given network design $(x)$ and disruption scenario $(d^s)$.  Since $(d^s)$  is a parameter, if an arc $(i,j)$ is disrupted in a given scenario $s$, we can replace $(x-d)$ with simply $0$.  In this scenario, the capacity of arc $(i,j)$ is zero regardless of whether the arc is constructed or not.  Alternatively, if arc $(i,j)$ is \emph{not} disrupted in scenario $s$,  we replace the right hand side of constraint (\ref{ef_arc_cap}) with simply $(x_{ij})$.  Finally, (\ref{ef_flow_nneg}--\ref{ef_x_integ}) are variable non-negativity and integrality constraints.

\textsc{ef} is a large-scale \textsc{milp}, with $O(|E|)$ binary variables and $O((|E|+|N|)\cdot S)$ continuous variables, and  $O((|E|+|N|)\cdot S)$ constraints.  Additionally, the linear programming relaxation of  (\ref{ef}) will likely be weak, as several  values of $(x)$ will be fractional at optimality for system components whose capacities are not fully utilized.  Therefore, solving (\ref{ef}) directly via a branch-and-bound approach will most likely be intractable for all but the smallest instances (which we demonstrate via computational experiments in section \ref{sec_comp_exp}). 

\subsection{A Benders Decomposition}

 The \textsc{ef}  formulation has a large number of variables and constraints.  For instances involving large networks and/or  a large disruption budget (and thus a large number of disruption scenarios $S$), formulation  (\ref{ef}) can become computationally intractable.  In this section, we present an alternative formulation with only $|E|$ binary variables $(x)$ and one continuous variable $\theta$ but possibly an extremely large number of constraints.  We use linear programming duality to generate valid inequalities for the projection of the natural formulation onto the space of the $(x)$ variables.  In essence, we use a variant of  Benders Decomposition in which we generate valid inequalities corresponding to ``optimality'' cuts.

Given a network design $(x)$ and disruption scenario $(d^s)$, consider the following linear program, denoted as the \emph{primal subproblem} \textsc{psp}$(x,d^s)$, to determine a flow $(f^s)$ that minimizes the demand shortage $\gamma^s$  under disruption scenario $(d^s)$.

\begin{subequations}\label{psp}
\begin{eqnarray}
 &\hskip -1.1cm   \theta^s(x,d^s) =  \ \displaystyle{\min_{f^s,\gamma^s}}  &  \gamma^s \label{psp_obj}\\
  & \hskip -1cm \mathrm{s.t.} \quad  (\alpha^s) &  \hskip -0.3cm  \sum_{j:(i,j)\in A}  f_{ij}^s - \sum_{j:(j,i)\in A } f_{ji}^s   = b_i (1-\gamma^s)  \nonumber \\ 
&&\hskip 4.1cm \forall i \in N \label{psp_flow_bal} \\
& \hskip -0.2cm(\beta^s)&  \hskip -0.3cm  f_{ij}^s, f_{ji}^s \leq u_{ij} (x_{ij} - d_{ij}^s) \  \forall \{i,j\} \in E \label{psp_arc_cap}\\
&& \hskip -0.3cm f_{ij}^s \geq 0 \quad \forall (i,j) \in A  \label{psp_flow_nneg}\\
&& \hskip -0.3cm \gamma^s \geq 0  \label{psp_gamma_nneg}
\end{eqnarray}
\end{subequations}

The objective (\ref{psp_obj}) is to minimize demand shortages given the prescribed network design $(x)$ and disruption scenario $s$.  Let the optimal objective value be given by $\theta^s(x,d^s)$.  Clearly, if $\theta^s(x,d^s)>0$, then there does not exist a feasible flow satisfying all demands, and if $\theta^s(x,d^s)=0$, then a feasible flow does exist satisfying all demand.  

In turn, we can formulate the dual of this problem as \textsc{dsp}$(x,d^s)$:

\begin{subequations}\label{dsp}
\begin{eqnarray}
&&\hskip -1.5cm \max_{\alpha^s, \beta^s} \quad \sum_{i \in N}  b_i\alpha_i^s + \sum_{\{i,j\} \in E} u_{ij}(x_{ij} - d_{ij}^s) (\beta_{ij}^s + \beta_{ji}^s)  \\
&\textmd{s.t.}&  \alpha_i^s - \alpha_j^s + \beta_{ij}^s \leq 0\quad \forall (i,j) \in A \label{dsp_f}\\
&&\sum_{i \in N}b_i \alpha_i^s \leq 1 \quad \forall i \in N  \label{dsp_gamma} \\
&&\beta_{ij}^s \leq 0 \quad \forall (i,j) \in A \label{dsp_beta_neg}
\end{eqnarray}
\end{subequations}

Since \textsc{psp}$(x,d^s)$ has a finite optimal solution (in the worst case, $\gamma^s$ is equal to 1 for all $s$) the dual \textsc{dsp}$(x,d^s)$ also has a finite optimal solution and, by strong duality, the optimal solutions coincide.  Therefore (\ref{dsp}) has a finite optimal solution, and the in fact, an optimal extreme point. Thus we can re-write the sub-problem as 

\begin{eqnarray}
 & \hskip -4cm \displaystyle \theta^s(x,d^s) = \max_{\ell = 1,\cdots, L} \Bigg(\sum_{i \in N}  b_i\alpha_i^\ell \nonumber \\
 &\hskip 2cm \displaystyle  +\!\!  \sum_{\{i,j\} \in E} \! u_{ij}(x_{ij}\! - \! d_{ij}^s)(\beta_{ij}^\ell + \beta_{ji}^\ell)\Bigg)
\end{eqnarray}

\noindent where $L$ is the set of extreme points corresponding to the polyhedron characterized by inequalities (\ref{dsp_f})--(\ref{dsp_beta_neg}). Note that the feasible region of $(\ref{dsp})$ does not depend on the network design $(x)$ or the disruption scenario $(d^s)$, which only affects the objective function.  

Alternatively, $\theta^s(x,d^s)$ is the smallest number $\theta^s$ such that
\begin{eqnarray}
&\displaystyle \sum_{i \in N}  b_i\alpha_i^\ell  +\!\!  \sum_{\{i,j\} \in E} \! u_{ij}(x_{ij}\! - \! d_{ij}^s) (\beta_{ij}^\ell+\beta_{ji}^\ell)  \leq \theta^s \nonumber \\
& \hskip 5cm \ell=1,\cdots,L \label{opt_cut}
\end{eqnarray}

Observing that $\theta^s \leq \theta$ for all $s$ and using (\ref{opt_cut}), we can reformulate (\ref{ef}) as

\begin{subequations}\label{proj_form}
\begin{eqnarray}
  &\displaystyle{\min_{x, \theta}}  & \sum_{\{i,j\} \in E} c_{ij}x_{ij} + p \hskip 0.02cm  \theta \label{rmp_obj}\\
&\mathrm{s.t.} &\sum_{i \in N}  b_i\alpha_i^\ell  +\!\!  \sum_{\{i,j\} \in E} \! u_{ij}(x_{ij}\! - \! d_{ij}^s) (\beta_{ij}^\ell + \beta_{ji}^\ell) \leq \theta \nonumber \\
&& \hskip 1cm \ell=1,\cdots,L, \ s=1,\cdots,S \label{rmp_opt_cut} \\  
&&x_{ij} \in\{0,1\} \quad \forall (\{i,j\} \in E \label{rmp_x_integ}\\
&&\theta \geq 0 \label{rmp_theta_nneg}
\end{eqnarray}
\end{subequations}

Formulation (\ref{proj_form}) has an exponential number of constraints, so we solve this via Benders decomposition (\textsc{bd}).  At a typical iteration of \textsc{bd}, we consider the \emph{relaxed master problem} (\textsc{rmp}), which has the same objective as problem  (\ref{proj_form}) but involves only a small subset of the constraints.  We briefly outline \textsc{bd} below.  For a detailed treatment of \textsc{bd} please refer to \cite{Benders1962}.

Let $t$ be the iteration counter and let the initial \textsc{rmp} be problem (\ref{proj_form}) without any constraints (\ref{rmp_opt_cut}).

\begin{algorithm}[H]
\caption{\emph{Benders Decomposition}}
\begin{algorithmic}[1]
\State $t\gets 0$
\State solve \textsc{rmp} and let ($x^t$, $\theta^t$) be the solution
\State \textbf{for} {$s=1,\cdots, S$}
\State \hskip 0.4cm \textbf{if} {\textsc{dsp}$(x^{t},d^s)> \theta^t$}
\State \hskip 0.8cm add cut $(\ref{opt_cut})$ to \textsc{rmp}
\State \hskip 0.4cm \textbf{end if}
\State \textbf{end for}
\State \textbf{if} {$\forall \ s=1,\cdots,S$, \ {\textsc{dsp}$(x^{t},d^s)\leq \theta^t$}}
\State \hskip 0.4cm $(x^t, \theta^t)$ is optimal (\textsc{exit})
\State \textbf{else}
\State \hskip 0.4cm $t\gets t+1$ and \textsc{goto} step 2
\State \textbf{end if}
\end{algorithmic}
\end{algorithm}

By using a Benders reformulation, we are able to decompose the extremely large \textsc{milp} (\ref{ef}) into a master problem and multiple subproblems, one for each disruption scenario.  This enables us to solve larger instances, which would not be possible by a direct solution of the \textsc{ef} formulation.  However, the extremely large number of disruption scenarios make direct application of Benders ineffective for instances with large networks and/or a large disruption budget.  In the next section, we develop a custom cut-generation based algorithm that evaluates all possible disruption scenarios implicitly using a separation oracle.

\subsection{Delayed Scenario Generation}

\begin{figure}[h]
   \centering
    \includegraphics[width=0.5\textwidth, angle=0]{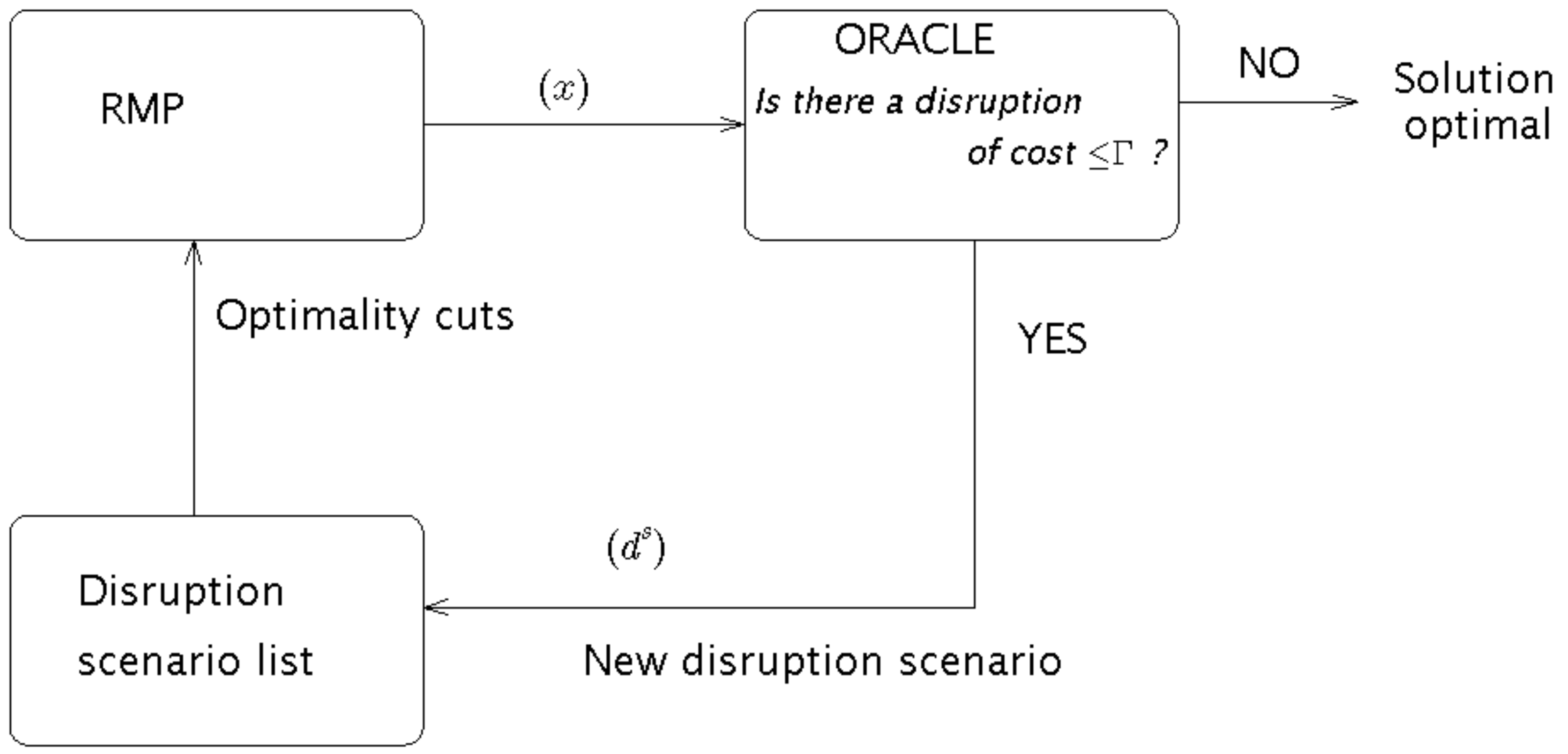}
    \caption{\footnotesize Proposed implicit optimization approach.} \label{dsg_flowchart}
\end{figure}

In practice, networks may be extremely large and  the number of disruption scenarios may be too large to be considered explicitly (even with a \textsc{bd} approach).  Our goal is to instead use an \emph{oracle} that implicitly evaluates all disruption scenarios and either identifies a violated one (with cost $\leq \Gamma$) or provides a certificate that no such disruption scenario exists.  If such a disruption exists, we add this disruption scenario to a disruption list and proceed to solve the updated list of scenario sub-problems.  If no such disruption exists, then the current network design $(x)$ is optimal and we terminate.  The proposed implicit optimization approach is summarized by the flowchart in Fig. \ref{dsg_flowchart}.

We next describe the oracle. Recall that we are trying to solve bi-level network disruption problem (\ref{ndp_bi-level}), i.e. given an $(x)$, we want the $(d)$ that maximizes the minimized shortage.  We state (without proof, for the sake of brevity) that this problem can be reformulated as (\ref{general_ndp}).

\begin{theorem}
\label{ndp_milp}
Problem (\ref{ndp_bi-level}) has an equivalent \textsc{milp} formulation as follows:
\begin{subequations}\label{general_ndp}
\begin{eqnarray}
& \hskip -1cm \textsc{ndp}(x) = \displaystyle \max_{d,\alpha,\beta} & \hskip -0.2cm \sum_{i \in N}  b_i\alpha_i + \sum_{\{i,j\} \in E} u_{ij}x_{ij} (\beta_{ij} + \beta_{ji}) \label{general_ndp_obj} \\
&\hskip -2.7cm \mathrm{s.t.}&\begin{array}{c}
    \hskip -1.8cm  \alpha_i - \alpha_j + \beta_{ij} -\hat{\beta}_{ij} d_{ij} \leq \\
    \hskip -1.8cm  \alpha_j - \alpha_i + \beta_{ji} -\hat{\beta}_{ji} d_{ij} \leq 
  \end{array}
  \left \}
    \begin{array}{c}
      0 \quad \forall \{i,j\} \in E
    \end{array}
  \right.\\
&&\hskip -1.65cm \sum_{i \in N}-b_i \alpha_i \leq 1 \quad \forall i \in N   \\
&&\hskip -1.65cm r^Td \leq \Gamma \label{general_ndp_budget}\\
&&\hskip -1.65cm \beta_{ij} \leq 0 \quad \forall (i,j) \in A \\
&&\hskip -1.65cm d_{ij} \in \{0,1\} \quad \forall \{i,j\} \in E
\end{eqnarray}
\end{subequations}

\end{theorem}

This is a standard \textsc{milp} with $O(|E|)$ binary variables, $O(|N|+|E|)$ continuous variables, and $O(|E|+|N|)$ constraints that can be solved much more easily than the exponential number of individual disruption scenarios. Similar reformulations for bi-level problems appear in \cite{Brown2006},\cite{Smith2007}. 
%


\begin{algorithm}[H]
\caption{\emph{Delayed Scenario Generation}}
\begin{algorithmic}[1]
\State $t\gets 0$
\State solve \textsc{rmp} and let ($x^t$, $\theta^t$) be the solution
\State solve \textsc{ndp}$(x^t)$ let $(d^t,\alpha^t, \beta^t)$ be the solution
\State \textbf{if} { {$ \textsc{ndp}(x^t)> \theta^t$} } \Comment{($d^t$) new scenario}
\State \hskip 0.4cm $t\gets t+1$ and add ($d^t$) to scenario list 
\State \hskip 0.4cm \textbf{for} {$s=1,\cdots, t$}
\State \hskip 0.8cm \textbf{if} { \textsc{dsp}$( x^{t-1},d^s) > \theta^t$}
\State \hskip 1.2cm add cut $(\ref{opt_cut})$ to \textsc{rmp}
\State \hskip 0.8cm \textbf{end if}
\State \hskip 0.4cm \textbf{end for}
\State \hskip 0.4cm \textsc{goto} step 2
\State \textbf{else}
\State \hskip 0.4cm $(x^t, \theta^t)$ is optimal (\textsc{exit})
\State \textbf{end if}
\end{algorithmic}
\end{algorithm}

At each iteration, either a new disruption is identified and added to the scenario list, leading to the generation of a new (violated) Benders cut, or the algorithm terminates with the current solution being optimal (if no new disruption is found).  Since the number of possible disruptions is finite, \textsc{dsg} will terminate after a finite number of iterations.  

Although  \textsc{dsg} is finitely convergent, the efficiency of the algorithm will largely be dictated by the strength of the \textsc{rmp} and \textsc{ndp} formulations and the order in which the disruption scenarios are identified.  In the next section, we exploit the problem structure to take advantage of special network properties. Specifically, we leverage the fact that max-flow equals min-cut -- to find a feasible disruption for which the min cut is below a certain threshold. If this min cut translates to a max flow less than the total demand of our network, then the network is infeasible -- i.e. we have identified a disruption scenario that leads to demand shortage. Using this insight, we derive a strong formulation of \textsc{ndp} and used this formulation to demonstrate the efficacy of \textsc{dsg} through extensive computational experiments in section \ref{sec_comp_exp}.


\subsection{A Strong \textsc{ndp} Formulation}

The performance and scalability of the  approach presented in the previous section directly depends on our ability to  identify   disruptions efficiently.  In this section, we present a strong formulation that  is designed  specifically for  inhibiting flows in networks.   This formulation has been previously used in~\cite{Pinar2010, Burch2003}. Here we  describe the formulation for completeness.  

The network disruption problem  seeks to  minimize the maximum flow of a network in a cost optimal way.  
Instead of  solving the bilevel optimization problem (\ref{ndp_bi-level}) or its corresponding \textsc{milp} reformulation (\ref{general_ndp}),  we use the duality between maximum flow and minimum cut, and  reformulate  the  problem to  remove lines from a given network so that there exists a cut  in the resulting network, whose capacity is less than the total demand in the system. This means that in the resulting system the demand cannot be satisfied completely and thus  we have a disruption within  the given budget.  
 
Let us add two special nodes to the network: a  source $s$ and a terminal  $t$, and connect each supply node $i$ to $s$ with a capacity equal to  the  supply capacity of the node $u_{si}=b_i$. Similarly, we connect  each demand node $j$ to the terminal $t$ so that $u_{jt}=b_j$.  We will refer to the new arc set that contains these new arcs as $E'$.  A cut in a network is defined by a bipartitioning of its nodes as  $N_1$ and $N_2=N\setminus N_1$ such  that $s\in N_1$ and $t\in N_2$.  The capacity of  this cut is defined by  the cumulative capacity  of arcs from $N_1$ to $N_2$, i.e.,   $\displaystyle \sum_{i\in N_1,j\in N_2}u_{ij}$.  
 We define a binary variable $\rho_i$ for each node $i \in N$,  so that
\[
\rho_i = \left\{ \begin{array}{ll}  0, & {\rm if}\;  i \in N_1, \\
                                                    1, & {\rm if}\;
                                                    i \in N_2,
 \end{array}\right.
\]
We also define a binary variable $\omega_{ij}\in [0,1]$ for each  arc. We will later show that $\omega_{ij}$  only takes values in $\{0,1\}$.
Binary variables $d_{ij}$ identify lines that are disrupted. We  can formulate the strong \textsc{ndp} as follows.

%
\begin{subequations}\label{mod_strong_ndp}
 \begin{align}
 \min_{\rho, \omega, d}\quad & \quad \sum u_{ij}w_{ij}  \label{mod_strong_ndp_obj}\\
  &  r^Td\leq \Gamma \label{eq:milp:cost} \\
  & \rho_i -\rho_j+ \omega_{ij}+d_{ij} \geq 0,  \label{eq:inhibit:cut2} \\
  & \rho_s=0,  \label{eq:inhibit:s} \\
  & \rho_t=1;  \label{eq:inhibit:t} \\
  & \rho_i\in \{0,1\} \;\;\; \forall  i \in N  \label{eq:inhibit:binary1} \\
  & d_{ij} \in \{ 0,1\}\;\;\;  \forall \{i,j\}\in E \label{eq:inhibit:binary2} \\
  & \omega_{ij} \in [0,1]\;\;\; \forall  \{i,j\}\in E' \label{eq:inhibit:binary3}
   \end{align}
   \end{subequations}
The objective function (\ref{mod_strong_ndp_obj}) measures the  capacity of the cut, and if it is below  the  total demand,  the $d$ variables identify a new disruption scenario that fails the network. Constraint (\ref{eq:milp:cost})  guarantees that the attack is  within budget.  Constraint (\ref{eq:inhibit:cut2}) ensures that  any arc directed from the source side $N_1$ to  the terminal side $N_2$ is either disrupted ($d_{ij}=1$) or  contributes to the cut  $(\omega_{ij}=1)$. Note that $\omega_{ij}=0$ in any other case due to the objective function, thus $\omega$ takes only binary values  even though it is a continuous variable.   
This formulation  offers significant improvements over the general \textsc{ndp} formulation (\ref{general_ndp}).  In our computational experiments, we observed formulation (\ref{mod_strong_ndp}) to run 3 to 5 times faster than the general \textsc{ndp} formulation, even for small instances, with the improvement gap increasing for larger problems.   The NDPs can be solved within seconds even for  very large problems and generous budgets. The reason for this speedup is that  the gap between this formulation and its  relaxation is very small, as has been analyzed in \cite{Burch2003}.

\section{Computational Experiments}\label{sec_comp_exp}

To test the performance of our proposed approach, we conducted computational experiments on a number of test cases, under a variety of parameters. All experiments were run on a machine with four quad-core 2.93G Xeon with 96G of memory.  For all computational experiments,  a single CPU and up to 8GB of RAM was allocated. CPLEX 11.2/Concert Technology v.27 was used for solving all mathematical programs.

Altogether, we considered twelve problem instances. We considered three different networks, derived from the \textsc{ieee} 30-node, 118-node, and 179-node test systems. For each of these networks, we considered four different disruption budgets. Specifically, we limited the attacker to disrupting at most one, two, three, or four arcs in the network. In all twelve instances, the power flow was approximated by a simple transportation model.

Table \ref{tab1} allows us to compare the run times for the three different approaches to solving the problem. For each of the twelve instances, $N$ provides the number of disruptable arcs.  For each test instance, we replicate the set of existing arcs multiple times to ensure that we have enough candidate arcs to yield a feasible network design.  Next, $k$ is the maximum number of arcs that can be disrupted by the attacker. In this case, the attack budget is given by $k$ and the cost of disruption is one for all arcs. The number of possible scenarios, i.e. the size of the set of feasible disruptions, appears in the next column. Note that for all but the smallest instances, this is a very large number -- approaching the hundreds of billions in the largest case.

 \begin{table}
\caption{Run times for different solution approaches}
\centering
\begin{tabular}{c c c c c c c}
\hline \hline
Test 		&		& 		&  No. poss. 	& 		& 		& 				\\ 
Systems 	&$N$ 	& $k$	&  scen. 	& \textsc{ef} 	& \textsc{bd}	& 	\textsc{dsg}			\\[0.5ex]
\hline
1	&	82		&	1	&	82 			& 0 			& 0 		& 0			\\
2	&	358		&	1	&	358 			& 20 			& 4 		& 4			\\
3	&	444		&	1	&	444 			& 33 			& 11 		& 19			\\
4	&	123		&	2	&	$>7K$ 		& $81722$ 	& 34 		& 1			\\
5	&	537		&	2	&	$>140K$ 		& x 			& $2142$ 	& 4			\\
6	&	666		&	2	&	$>200K$ 		& x 			& $5924$ 	& 174		\\
7	&	164		&	3	&	$>700K$ 		& x 			& $3045$ 	& 9			\\
8	&	716		&	3	&	$>60M$		& x 			& x 		& 398		\\
9	&	888		&	3	&	$>116M$ 		& x 			& x 		& 653		\\
10	&	205		&	4	&	$>72M$ 		& x 			& x 		& 67			\\
11	&	895		&	4	&	$>26B$ 		& x 			& x 		& 2708		\\
12	&	$1110$	&	4	&	$>63B$ 		& x 			& x		& $11999$	\\ [0.5ex]
\hline
\end{tabular}
\label{tab1}
\end{table}

The remainder of this table provides the run time (in seconds) for each instance under the three different approaches. Note that the first approach, the extensive form (\textsc{ef}), can only solve the smallest of instances. This is because of the sheer size of the problem Ð for each contingency, a full network flow problem must be embedded in the formulation. As the number of contingencies grows, this quickly becomes intractable.

The second approach, \textsc{bd}, bypasses this problem via a Benders decomposition, with corresponding delayed cut generation. However, this still suffers from the exponential growth in the number of disrupton scenarios -- for each contingency, a dual subproblem must be solved to check for violated Benders cuts to add to the master problem. We see that larger problem instances can be solved, relative to \textsc{ef}, but that the \textsc{bd} approach nonetheless cannot solve the largest problem instances.

In the \textsc{dsg} approach (using the strong \textsc{ndp} formulation (\ref{mod_strong_ndp})), we see that all instances of the problem can be solved, in almost all cases in under an hour and frequently in only a few minutes. This is a result of the combination of the strength of the Benders cuts, enabling the problem to be solved in a very limited number of iterations, and also the fact that we are able to implicitly evaluate the contingencies in order to identify a violated contingency and then quickly find its corresponding Benders cut.

Table \ref{tab2} provides us with further evidence to support this. For each instance, we see the total number of possible disruption scenarios and then the number of disruption scenarios for which corresponding cuts were actually generated (this is the total number of iterations). Clearly, it is a very tiny fraction of the possible number of disruptions, which is critical to the tractability of the approach. The remaining columns of this table breakdown the total run time by time spent on the three components of the algorithm -- the restricted master problem (\textsc{rmp}), which identifies a candidate network design $(x)$; the network disruption problem (\textsc{ndp}), which identifies a contingency that cannot be overcome by the current network design; and the scenario subproblems (\textsc{sp}), which generates the Benders cuts.  It is interesting to note that no one category consistently dominates the total time.

\begin{table}
\caption{\textsc{dsg} runtime breakdown}
\centering
\begin{tabular}{c c c c c c c c c}
\hline \hline
Test 		&  No. poss.	&	No. eval.		& Total		& \textsc{rmp}	& \textsc{ndp}	&\textsc{sp}			\\ 
Systems 	&  scen. 		& 	 scen.		& time 		& time	& time	&time		\\[0.5ex]
\hline
1		&82 			& 	3 			& 0			&0		&0		&0			\\
2		&358 		& 	17 			& 4			&0		&2		&1			\\
3		&444 		& 	51 			& 19			&1		&7		&10			\\
4		&$>7K$ 		& 	15 			& 1			&0		&1		&0			\\
5		&$>140K$ 	& 	58 			& 4			&3		&26		&12			\\
6		&$>200K$ 	& 	158 			& 174		&6		&50		&118		\\
7		&$>700K$ 	& 	43 			& 9			&2		&5		&2			\\
8		&$>60M$		&	128 			& 398		&25		&303	&70			\\
9		&$>116M$ 	& 	284 			& 653		&21		&193	&439		\\
10		&$>72M$ 		& 	156 			& 67			&7		&23		&37			\\
11		&$>26B$ 		& 	359 			& 2708		&$399$	&1698	&612		\\
12		&$>63B$ 		& 	899 			& $11999$	&$4939$	&$1822$	&$5237$		\\ [0.5ex]
\hline
\end{tabular}
\label{tab2}
\end{table}

\begin{figure}[h]
   \centering
   \vskip -0.75cm
    \includegraphics[width=0.5\textwidth, angle=0]{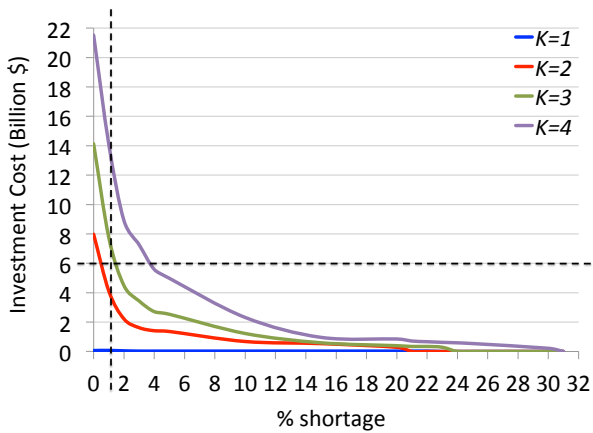}
    \vskip -1.0cm
    \caption{\footnotesize Investment cost and demand shortage trade-off}
    \label{fig1}
\end{figure}

Finally, the fact that our proposed approach enables us to find \emph{optimal} solutions to the twelve problem instances in tolerable run time provides us with the ability to also conduct analysis on the quality of the solutions.  Figure \ref{fig1} provides some introductory insights into the trade-off between investment in security and potential shortfalls for the \textsc{ieee} 118-node system.  The figure shows the total investment cost (to secure the network) as a function of allowable demand shortages for each of the four disruption budgets.  Observe, that the \textsc{ieee} 118-node system is largely secure for $k=1$.  For higher disruptions, a significant cost is incurred for completely securing the system (zero shortage).  The vertical dashed line indicates that significant investment cost savings are attainable if a one percent demand shortage is permitted. On the other hand, if we restrict the network investment budget to \$6 billion (highlighted by the dashed horizontal line), we can assess the security of the network against the various disruption budgets.

\section{Conclusion}\label{sec5}
Network security is critical to society in many forms -- telecommunications, transportation, water, and electricity, to name a few. Concerns exist over both malicious, intentional attacks and random failures tied to systems that are rapidly approaching the boundaries of their capacity. With the increased size and complexity of these networks, direct human analysis and operational control are no longer possible. Instead, automated tools are required to ensure the security of these networks.

Substantial research has been invested in identifying network vulnerabilities, often through bi-level optimization approaches. In this paper, we have expanded this research with the goal of identifying cost-effective ways to secure a network against such vulnerabilities. We first pose this problem as an explicit tri-level optimization problem. This allows us to clearly state the problem, but is not tractable in practice because of the enormous size -- the formulation must simultaneously capture separate network flow problems for each of an exponential number of disruption scenarios. We next consider a Benders decomposition-based approach to the problem. This is certainly an improvement, with Benders cuts replacing the exponential number of embedded network flow problems. In order to identify the Benders cuts, however, each of the exponential number of contingencies must be evaluated, which still greatly limits the size of the problem instance that can be solved in an acceptable run time.

To overcome this problem, we propose the new approach of using a separation problem to implicitly consider the full set of contingencies in a much more efficient fashion. Each iteration of the separation problem yields a contingency that, relative to the current network design, would lead to infeasible network flows. For such a contingency, we can then directly generate the associated Benders cut.  This new approach is demonstrated to be much more effective in solving problem instances of a realistic size in acceptable run times.

This research also lays the foundation for future, more general research on tri-level optimization problems and, in particular, survivable network design. In particular, we are focused on first expanding the simple network flow problem in the third stage to more complex (and realistic) power flow models. Beyond this, we will be working to generalize the approach to the broad class of problems where the third stage problem is not restricted only to be a linear program. Finally, we will consider the case where third stage decisions have non-linear form. 

\section*{Acknowlegments}
This work was funded by the applied mathematics program at the United States Department of Energy and Laboratory Directed Research and Development Program of Sandia National Laboratories, a multiprogram laboratory operated by Sandia Corporation, a wholly owned subsidiary of Lockheed Martin Corporation, for the United States Department of Energy's National Nuclear Security Administration under  contract DE-AC04-94AL85000. We thank Jean-Paul Watson for many productive discussions.

\end{document}